\documentclass[a4paper,11pt]{article}

\usepackage{amsmath}
\usepackage{amssymb}  
\usepackage{dsfont} 
\usepackage{vmargin}

\newcommand{\R}[1]{r_{\epsilon}^{#1}} 
\newcommand{\y}[1]{y_{\epsilon}^{#1}} 

\newcommand{\E}{\mathbb{E}} 
\newcommand{\p}{\mathbb{P}} 
\newcommand{\vc}[1]{\mathcal{V}_{#1}}  
\newcommand{\et}[1]{\tilde{\mathbb{E}}_{(x,#1,0)}}
\newcommand{\T}[1]{\tilde{#1}}
\newcommand{\id}[1]{\mathds{1}_{\{#1\}}} 
\newcommand{\e}{\mathrm{e}} 
\DeclareMathOperator*{\argmin}{arg\,min} 

\begin{document}

\newtheorem{thm}{Theorem}[section]
\newtheorem{prop}[thm]{Proposition}
\newtheorem{lem}[thm]{Lemma}
\newtheorem{corollary}[thm]{Corollary}
\newtheorem{definition}[thm]{Definition}
\newtheorem{rem}[thm]{Remark}
\newtheorem{conjecture}[thm]{Conjecture}
\newtheorem{assum}[thm]{Assumption}
\newtheorem{condition}[thm]{Condition}

\title{Optimal strategies for impulse control of piecewise deterministic Markov processes\thanks{This work was supported by the French National Agency of Research (ANR), grant
PIECE, number ANR-12-JS01-0006}} 

\author{ \mbox{ }
\\ \\ B. de Saporta
\\ 
\small Universit\'e de Montpellier, France \\
\small IMAG, CNRS UMR 5149, France \\
\small INRIA Bordeaux Sud Ouest\\
\small e-mail :  benoite.de-saporta@umontpellier.fr
\\
\and
\\ \\ F. Dufour
\\ 
\small Bordeaux INP, France\\
\small IMB, CNRS UMR 5251, France \\
\small INRIA Bordeaux Sud Ouest \\
\small e-mail : dufour@math.u-bordeaux1.fr 
\\
\and
\\ \\ A. Geeraert
\\ 
\small INRIA Bordeaux Sud Ouest\\
\small Universit\'e de Bordeaux, France\\
\small IMB, CNRS UMR 5251, France\\
\small e-mail : alizee.geeraert@inria.fr
}

\maketitle
\begin{abstract}     
This paper deals with the general discounted impulse control problem of a piecewise deterministic Markov process. We investigate a new family of $\epsilon$-optimal strategies.
The construction of such strategies is explicit and only necessitates the previous knowledge of the cost of the no-impulse strategy. 
In particular, it does not require the resolution of auxiliary optimal stopping problem or the computation of the value function at each point of the state space.
This approach is based on the iteration of a single-jump-or-intervention operator associated to the piecewise deterministic Markov process.  

\end{abstract}



\section{Introduction}

The aim of this paper is to propose a new family of $\epsilon$-optimal strategies for the impulse control problem of piecewise deterministic Markov processes (PDMPs) defined by O.L.V. Costa and M.H.A. Davis in \cite{Costa89}. We consider the infinite horizon expected discounted impulse control problem where the controller instantaneously moves the process to a new point of the state space at some controller specified time.

Piecewise deterministic Markov processes have been introduced by M.H.A. Davis in \cite{Davis84,Davis93} as a general class of stochastic hybrid models. These processes have two variables: a continuous one representing the physical parameters of the system and a discrete one which characterizes the regime of operation of the physical system and/or the environnement. The process depends on three local characteristics: the flow, the jump intensity and the Markov kernel. 
The path of a PDMP consists of deterministic trajectories punctuated by random jumps. 
 Roughly speaking, the behavior of the process is the following. Starting from a point of the state space, the PDMP follows a deterministic trajectory determined by the flow, until the first jump time. 
 This time is drawn either in a Poisson like fashion following the jump intensity or deterministically when the process hits the boundary of the state space. 
 The new position and regime of the PDMP after the jump are selected by the Markov kernel.
Then the process follows again a deterministic trajectory until the next jump time and so on. 
There are many and diverse applications of PDMPs for example in queuing or inventory systems, insurance, finance, maintenance models  \cite{bauerle11,Dassios89,Davis93} or in data transmission \cite{Chafai10} and in biology \cite{doumic2015,pakdaman10}. The interested reader can also refer to \cite{dd15} for some applications in reliability.

Impulse control corresponds to the following situation: the process runs until a controller decides to intervene on the process by instantaneously moving the process to some new point of the state space. Then, restarting at this new point, the process runs until the next intervention and so on. 
Many authors have considered impulse control for PDMPs, either by variational inequality \cite{Lenhart89,Gatarek92,Dempster95} or by value improvement \cite{Costa89}. 
The simplest form of impulse control is optimal stopping, where
the decision maker selects only one intervention time when the process is stopped.
Optimal stopping for PDMPs has been studied in particular in \cite{Gugerli86,Costa88,BF2010}. 

For a general infinite horizon expected discounted impulse control problem, a strategy consists in two sequences of controller-specified random variables defining the intervention times and new starting points of the process.
Solving this problem involves finding a strategy that minimizes the expected sum of discounted running and intervention costs up to infinity.
The minimal cost is called the value function.
In general, optimal strategies do not exist. Instead we consider $\epsilon$-optimal strategies, i.e. strategies which cost differs from the value function at most of $\epsilon$. 


There exists an extensive literature related to the study of the optimality equation associated to expected discounted control problems but few works are devoted to the characterization of $\epsilon$-optimal strategies. The objective of the paper is to explicitly construct such strategies. An attempt in this direction has been proposed by O.L.V. Costa and M.H.A. Davis in \cite[section 3.3]{Costa89}. Roughly speaking, one step of their approach consists in solving an optimal stopping problem which makes this technique quite difficult to implement. Furthermore the knowledge of the optimal value function is required.

We propose a construction of an $\epsilon$-optimal strategy which necessitates only the knowledge of the cost of the non-impulse strategy and without solving technical problems preliminary. This construction is based on the iteration of a single-jump-or-intervention operator associated to the PDMP. It builds on the explicit construction of $\epsilon$-optimal stopping times developed by U.S. Gugerli \cite{Gugerli86} for the optimal stopping problem. However, for the general impulse control problem, one must also optimally choose the new starting points of the process, which is a significant source of additional difficulties. It is important to emphasize that our method has the advantage of being constructive with regard to other works in the literature on impulse control problem.

This work is also the first step towards a tractable numerical approximation of $\epsilon$-optimal strategies. A numerical method to compute the value function is proposed in \cite{BF2012}. It is based on the quantization of an underlying discrete-time Markov chain related to the continuous process and path-adapted time discretization grids. Discretization of $\epsilon$-optimal strategies will be the object of a future work.

The paper is organized as follow. In Section~\ref{Sec2} we recall in section the definition of a PDMP and state the impulse control problem under study. In section \ref{sect_approximation}, we construct a sequence of approximate value functions. In section \ref{Sec4}, we build aa auxilliary process corresponding to an explicit family of strategies and we show in section that the cost of the controlled trajectories corresponds to the approximate value function built in section \ref{sect_approximation}. Technical details are gathered in the Appendix.

\section{Impulse control problem of PDMP} \label{Sec2}

We introduce first some standard notation before giving a precise definition of a piecewise deterministic Markov processes (PDMP) and of the impulse control problem we are interested in.

For $a,b\in\mathbb{R}$, $a\wedge b=\min(a,b)$ 
is the minimum 
of $a$ and $b$. By convention, set $\inf \varnothing =\infty$.
Let $\mathcal{X}$ be a metric space with distance $d_{\mathcal{X}}$. For a subset $A$ of $\mathcal{X}$, $\partial A$ is the boundary of $A$ and $\bar{A}$ its closure.
We denote $\mathfrak{B}(\mathcal{X})$ the Borel $\sigma$-field of $\mathcal{X}$ and $\textbf{B}(\mathcal{X})$ the set of real-valued, bounded and measurable functions defined on $\mathcal{X}$. 
For any function $w \in \textbf{B}(\mathcal{X})$, we write $C_w$ for its upper bound, that is
$C_w = \sup_{x\in\mathcal{X}} |w(x)|$.
For a Markov kernel $P$ on $(\mathcal{X},\mathfrak{B}(\mathcal{X}))$ and functions $w$ and $w'$ in $\textbf{B}(\mathcal{X})$, for any $x\in\mathcal{X}$, set $Pw(x)=\int_{\mathcal{X}} w(y) P(x,dy)$.

\subsection{Definition of PDMP}\label{def_pdmp}

Let $M$ be the finite set of the possible regimes or modes of the system. For all modes $m$ in $M$, let $E_m$ be an open subset of $\mathbb{R}^d$ endowed with the usual Euclidean norm $|\cdot|$. Set $E=\{(m,\zeta), m\in M, \zeta \in E_m\}$.
Define on $E$ the following distance, for $x=(m,\zeta)$ and $x'=(m',\zeta')\in E$,
\[|x-x'| = |\zeta-\zeta'| \id{m=m'}+\infty \id{m \neq m'}.\]
A piecewise deterministic Markov process (PDMP) on the state space $E$ is determined by three local characteristics:
\begin{itemize}
\item the flow $\Phi(x,t)=(m,\Phi_m(\zeta,t))$ for all $x=(m,\zeta)$ in $E$ and $t\geq 0$, where $\Phi_m: \mathbb{R}^d \times \mathbb{R}^+ \rightarrow \mathbb{R}^d$ is continuous such that $\Phi_m(\cdot,t+s)=\Phi_m(\Phi_m(\cdot,t),s)$, for all $t,s\in \mathbb{R}^+$. It describes the deterministic trajectory between jumps. 
We set $t^*(x)$ the time the flow takes to reach the boundary of $E$ when it starts from $x=(m,\zeta)$:
\[t^*(x)=\inf\{t>0: \Phi_m(\zeta,t)\in \partial E_m \}.\]
\item the jump intensity $\lambda:\bar{E}\rightarrow \mathbb{R}^+$ is a measurable function and has the following integrability property: for any $x=(m,\zeta)$ in $E$, there exists $\epsilon >0$ such that
\[ \int_0^{\epsilon} \lambda(m,\Phi_m(\zeta,t))dt < +\infty. \]
For all $x=(m,\zeta)$ in $E$ and $t\in[0;t^*(x))$, we set
\begin{equation}\label{def_Lambda}
\Lambda(m,\zeta,t)=\int_0^{t} \lambda(m,\Phi_m(\zeta,s)) ds.
\end{equation}
\item the Markov kernel $Q$ on $(\bar{E},\mathfrak{B}(\bar{E}))$ represents the transition measure of the process and allows to select the new location after each jump. It satisfies for all $x \in \bar{E}$, $Q(x,\{x\}\cup \partial E)=0$.
That means each jump is made in $E$ and changes the location and/or the mode of the process.
\end{itemize}
From these characteristics, it can be shown \cite[section 25] {Davis93} that there exists a filtered probability space $(\Omega,\mathcal{F}$, $\{\mathcal{F}_t\},\{\p_x\}_{x\in E})$ on which a process $\{X_t\}$ can be defined as a strong Markov process. The process $\{X_t\}$ has two components $X_t=(m_t,\zeta_t)$ where the first component $m_t$ is usually called the mode or the regime and the second component $\zeta_t$ is the so-called Euclidean variable. The motion of this process can be then defined iteratively as follows. 

Starting at an initial point $X_0=(m_0,\zeta_0)$ with $m_0\in M$ and $\zeta_0 \in E_{m_0}$, the first jump time $T_1$ is determined by the following survival equation
\begin{equation}\label{loi_T}
\p_{(m_0,\zeta_0)}(\{T_1>t\})=  \e^{-\Lambda(m_0,\zeta_0,t)} \id{t<t^*(m_0,\zeta_0)}.
\end{equation}
On the interval $[0,T_1)$, the process $\{X_t\}$ follows the deterministic trajectory $\zeta_t=\Phi_{m_0}(\zeta_0,t)$ and the regime $m_t$ is constant and equal to $m_0$. At the random time $T_1$, a jump occurs. A jump can produce either a discontinuity in the Euclidean variable $\zeta_t$ and/or change of mode. The process restarts at a new mode and/or position $X_{T_1}=(m_{T_1},\zeta_{T_1})$, according to the distribution $Q((m_0,\Phi_{m_0}(\zeta_0,T_1)),\cdot)$. An inter jump time $T_2-T_1$ is then selected in a similar way to equation \eqref{loi_T}, and on the interval $[T_1,T_2)$, the process follows the path $m_t=m_{T_1}$ and $\zeta_t=\Phi_{m_{T_1}}(\zeta_{Z_1},t-T_1)$. The process $\{X_t\}$ thus defines a PDMP on the state space $E$.

In order to avoid any technical problems due to the possible explosion of the process, we make the following standard assumptions \cite[section 24]{Davis93}, \cite[section 1.4]{dd15}.
\begin{assum} The mean number of jumps before an instant $t\in\mathbb{R}^+$ is finite, whatever the initial position of the process: for all $x\in E$ and $t \in \mathbb{R}^+$, $\E_x\Bigl[\sum_{n=1}^{\infty} \id{T_n \leq t}\Bigr]< \infty$.
\end{assum}
This first assumption implies in particular that $T_k \rightarrow \infty$ almost surely, when $k\rightarrow \infty$.
\begin{assum}\label{hyp2}
The exit time $t^*$ is Lipschitz-continuous and bounded by $C_{t^*}$.
\end{assum}
In most practical applications, the physical properties of the system ensure that either $t^*$ is bounded or the problem has a natural deterministic time horizon $T$. In the latter case, there is no loss of generality in considering that $t^*$ is bounded by this deterministic time horizon. This leads to replacing $C_{t^*}$ by $T$.

We need also Lipschitz assumptions on the jump rate $\lambda$ and the Markov kernel $Q$.
\begin{assum}\label{hyp3}
The jump rate $\lambda$ is in $\textbf{B}(E)$. It is bounded and there exists $[\lambda]_1\in\mathbb{R}^+$ such that for any $(x,y)\in E^2$, $u\in[0,t^*(x)\wedge t^*(y))$,
\[ |\lambda(\Phi(x,u))-\lambda(\Phi(y,u))| \leq [\lambda]_1 |x-y|.\]
\end{assum}
We define $\textbf{L}_{\Phi}(E)$ as the set of functions $w\in \textbf{B}(E)$ that are Lipschitz continuous along the flow i.e. the real-valued, bounded, measurable functions defined on $E$ and satisfying the following conditions:
\begin{enumerate}
\item For all $x\in E$, the map $w(\Phi(x,\cdot)):[0;t^*(x)) \rightarrow \mathbb{R}$ is continuous, the limit $\displaystyle\lim_{t \rightarrow t^*(x)} w(\Phi(x,t))$ exists and is denoted by $w(\Phi(x,t^*(x)))$,
\item there exists $[w]_1 \in \mathbb{R}^+$ such that for any $(x,y) \in E^2$, $t\in [0;t^*(x)\wedge t^*(y)]$, one has
\[|w(\Phi(x,t))-w(\Phi(y,t))| \leq [w]_1 |x-y|, \]
\item there exists $[w]_2 \in \mathbb{R}^+$ such that for any $x \in E$, $(t,s) \in [0;t^*(x)]^2$, one has
\[|w(\Phi(x,t))-w(\Phi(x,s))| \leq [w]_2 |t-s|, \]
\item there exists $[w]_* \in \mathbb{R}^+$ such that for any $(x,y) \in E^2$, one has
\[|w(\Phi(x,t^*(x)))-w(\Phi(y,t^*(y)))| \leq [w]_* |x-y|. \]
\end{enumerate}

\begin{assum}\label{hyp4}
The Markov kernel $Q$ is Lipschitz in the following sense: there exists $[Q]\in\mathbb{R}^+$ such that for any function $w\in \textbf{L}_{\Phi}(E)$ the Markov kernel $Q$ satisfies:
\begin{enumerate}
\item for any $(x,y)\in E^2$, and $t\in[0,t^*(x)\wedge t^*(y))$, we have
\[ |Qw(\Phi(x,t))-Qw(\Phi(y,t))| \leq [Q][w]_1|x-y|, \]
\item for any $(x,y)\in E^2$, we have
\[ |Qw(\Phi(x,t^*(x)))-Qw(\Phi(y,t^*(y)))| \leq [Q][w]_*|x-y|. \]
\end{enumerate}
\end{assum}
For notational convenience, we set $T_0=0$. The sequence $(Z_n)_{n\in\mathbb{N}}$, with $Z_n=X_{T_n}$, describes the post jump locations of the process $\{X_t\}$. The sequence $(S_n)_{n\in\mathbb{N}}$, with $S_n=T_n-T_{n-1}$ for $n\geq 1$ and $S_0=0$, gives the sojourn times between two consecutive jumps. We can see that the process defined by $\{\Theta_n\}$, where $\Theta_n=(Z_n,S_n)$ for all $n\in\mathbb{N}$, is a discrete Markov chain and it is the only source of randomness of the process $\{X_t\}$.

\subsection{Impulse control problem}\label{sect_impulse_control}

The formal probabilistic apparatus necessary to precisely define the impulse control problem is rather cumbersome, and will not be used in the sequel. Therefore, for the sake of simplicity, we only present a rough description of the problem. The interested reader is referred to \cite{Costa89} or \cite[section 54]{Davis93} for a rigorous definition.

A strategy $\mathcal{S}=(\tau_n,R_n)_{n\geq1}$ is a sequence of non-anticipative intervention times $(\tau_n)_{n\geq1}$ and non-anticipative $E$-valued random variables $(R_n)_{n\geq1}$ on a measurable space $(\bar{\Omega},\bar{\mathcal{F}})$. Between the intervention times $\tau_i$ and $\tau_{i+1}$, the motion of the process is determined by the characteristics of the PDMP $\{ X_t\}$ starting from $R_i$. If an intervention takes place at $x\in E$, the set of admissible points where the decision-maker can send the system to is denoted by $\mathbb{U} \subset E$. We suppose that the control set $\mathbb{U}$ is finite and does not depend on $x$. The cardinal of the set $\mathbb{U}$ is denoted by $u$:
$\mathbb{U}=\{y^i: 1 \leq i \leq u \}$.
The strategy $\mathcal{S}$ induces a family of probability measures $\p^{\mathcal{S}}_x$, $x\in E$, on $(\bar{\Omega},\bar{\mathcal{F}})$.
We consider that the strategy is admissible if it satisfies the conditions given in \cite[section 2.3]{Costa89}. We denote $\mathbb{S}$ the class of admissible strategies.

Associated to the strategy $\mathcal{S}\in\mathbb{S}$, we define the following discounted cost for a process starting at $x\in E$
\begin{equation}\label{cout_strategie}
\mathcal{J}^{\mathcal{S}}(x)=\E_x^{\mathcal{S}} \Bigl[ \int_0^{\infty} \e^{-\alpha s}f(\tilde{X}_s)ds + \sum_{i=1}^{\infty} \e^{-\alpha \tau_i} c(\tilde{X}_{\tau_i-},\tilde{X}_{\tau_i})\Bigr],
\end{equation}
where $\E_x^{\mathcal{S}}$ is the expectation with respect to $\p^{\mathcal{S}}_x$, $\{\tilde{X}_t\}$ is the process with interventions (for its construction, see \cite[section 2.2]{Costa89}) and $\alpha$ is a positive discount factor. The function $f$ corresponds to the running cost. It is a non negative function in $\textbf{L}_{\Phi}(E)$. The function $c$ is a continuous function on $\bar{E}\times\mathbb{U}$, where $c(x,y)$ corresponds to the intervention cost of moving the process from $x$ to $y$. We add some assumptions on the intervention cost $c$. 
\begin{enumerate}
\item There exists $[c]_1\in\mathbb{R}^+$ such that for any $(x,y)\in E^2$ and $u\in[0,t^*(x)\wedge t^*(y)]$ we have
\[\max_{z\in\mathbb{U}} |c(\Phi(x,u),z)-c(\Phi(y,u),z)| \leq [c]_1 |x-y|. \]
\item There exists $[c]_2\in\mathbb{R}^+$ such that for any $x\in E$ and $(t,s)\in[0,t^*(x)]^2$ we have
\[\max_{z\in\mathbb{U}} |c(\Phi(x,t),z)-c(\Phi(x,s),z)| \leq [c]_2 |t-s|. \]
\item There exists $[c]_*\in\mathbb{R}^+$ such that for any $(x,y)\in E^2$ we have
\[\max_{z\in\mathbb{U}} |c(\Phi(x,t^*(x)),z)-c(\Phi(y,t^*(y)),z)| \leq [c]_* |x-y|. \]
\item For any $(x,y)\in\bar{E}\times\mathbb{U}$, there exist $c_0,C_c\in\mathbb{R}^+_*$ such that
\[0<c_0\leq c(x,y)\leq C_c.\]
\item For any $(x,y,z)\in\bar{E}\times\mathbb{U}\times\mathbb{U}$, we have
\[ c(x,y)+c(y,z) \geq c(x,z).\]
It means that the cost of one single intervention will be less than or equal to the cost of two simultaneous interventions.
\end{enumerate}


The value function of the discounted infinite horizon impulse control problem is defined for all $x$ in $E$ by
\[\mathcal{V}(x)=\inf_{\mathcal{S}\in \mathbb{S}} \mathcal{J}^{\mathcal{S}}(x).\]
Associated to this impulse control problem, we define the following operators \cite{Costa89,BF2010}.
For $x\in E$, $t\geq 0$, $(v,w)\in\textbf{L}_{\Phi}(E)^2$, set
\begin{equation} \label{def_F}
\begin{split}
F(x,t)=& \int_0^{t \wedge t^*(x)} \e^{-\alpha s-\Lambda(x,s)} f(\Phi(x,s)) ds \\
=&\E_x \Bigl[\int_0^{T_1 \wedge t} \e^{-\alpha s} f(\Phi(x,s)) ds\Bigr],
\end{split}
\end{equation}
\begin{equation} \label{def_J}
\begin{split}
J(v,w)(x,t)=& \int_0^{t \wedge t^*(x)} \e^{-\alpha s-\Lambda(x,s)} \\
&\times\Bigl[f(\Phi(x,s))
+\lambda Qw(\Phi(x,s))\Bigr] ds \\
&+ \e^{-\alpha(t\wedge t^*(x))-\Lambda(x,t\wedge t^*(x))} v(\Phi(x,t\wedge t^*(x))) \\
=& F(x,t)+ \E_x \left[\right.\e^{-\alpha (t\wedge t^*(x))}  \\
&\times v(\Phi(x,t\wedge t^*(x))) \id{S_1 \geq t\wedge t^*(x)} \\
& +\e^{-\alpha S_1} w(Z_1)\id{S_1 < t\wedge t^*(x)}  \left.\right],
\end{split}
\end{equation}
\begin{equation} \label{def_K}
\begin{split}
Kw(x)=& \int_0^{t^*(x)} \e^{-\alpha s-\Lambda(x,s)} \Bigl[\Bigr. f(\Phi(x,s)) \\
&+\lambda Qw(\Phi(x,s))\Bigl.\Bigr] ds \\
&+ \e^{-\alpha t^*(x)-\Lambda(x,t^*(x))} Qw(\Phi(x,t^*(x))) \\
=&F(x,t^*(x))+ \E_x \left[\e^{-\alpha S_1} w(Z_1)\right].
\end{split}
\end{equation}
For $(v,w)\in\textbf{L}_{\Phi}(E)^2$, $\varphi$ defined on $\mathbb{U}$ and $x\in E$, set
\begin{equation} \label{def_M}
M\varphi(x)= \inf_{y\in \mathbb{U}} \{c(x,y)+\varphi(y)\},
\end{equation}
\begin{equation} \label{def_L}
L(v,w)(x)=\inf_{t \in \mathbb{R}^+} J(v,w)(x,t) \wedge Kw(x),
\end{equation}
\begin{equation} \label{def_Lronde}
\mathcal{L}w(x)=L(Mw,w)(x).
\end{equation}


As explained in \cite{Costa89}, operator $\mathcal{L}$ applied to $w$ is the value function of the single jump-or-intervention problem, with cost function $w$, and the value function $\mathcal{V}$ can be computed by iterating $\mathcal{L}$. More precisely, let $h$ be the cost associated to the no-impulse strategy:
\begin{equation}\label{def_h}
\forall x \in E, \qquad h(x)=\E_x\Bigl[ \int_0^{\infty} \e^{-\alpha s} f(X_s)ds \Bigr]. 
\end{equation}


Then we can recall Proposition 4 of \cite{Costa89}.
\begin{prop}\label{vrond}
Define $\mathcal{V}_0=h$ and $\mathcal{V}_{n+1}=\mathcal{L}(\mathcal{V}_{n})$ for all $n \geq 0$. Then for all $x$ in $E$, one has
\[\mathcal{V}(x)=\lim_{n \rightarrow + \infty} \mathcal{V}_{n}(x).\]
\end{prop}

Function $\mathcal{V}_{n}$ corresponds to the value function of the impulse control problem where at most $n$ interventions are allowed, and the process runs uncontrolled after the $n$-th jump-or-intervention.


The aim of this paper is to propose an explicit $\epsilon$-optimal strategy for the impulse control problem described above.
Our approach is based on the sequence of value functions $(\mathcal{V}_{n})$. The first step, detailed in Section \ref{sect_approximation}, is to iteratively construct an approximation of the value functions $(\mathcal{V}_{n})$ based on an approximation of operator $\mathcal{L}$. Roughly speaking, replace the infimum of $J$ by a point at which the value of operator $J$ is at distance less than $\epsilon$ from the infimum of $J$. We prove that this sequence of approximate value function still converges to $\mathcal{V}$. The second step, given in Section \ref{sect_def_strategie}, consists in defining an auxiliary PDMP that can be interpreted as a controlled version of the original PDMP for a given explicit strategy. Then we establish in Section \ref{sect_optimal} that the cost of this strategy is exactly the approximate value function of Section \ref{sect_approximation}. Hence, we obtain an explicit $\epsilon$-optimal strategy.

\section{Approximate value fonctions }\label{sect_approximation}

We now construct an approximation of the sequence $(\mathcal{V}_{n})$.
%
%
When iterating operator $\mathcal{L}$, the controller selects the option generating the lowest cost between waiting for the next jump time (operator $K$ wins the minimization) and making an intervention (operator $J$ wins the minimization). In the latter case, the best time to intervene is given by the time achieving the infimum in $\inf J$. Hence, to construct an approximation of the value functions, one may use the above rules with an approximation of the infimum in $\inf J$.


Set $V_0=h$ and $\epsilon>0$. First, we define for any $x$ in $E$:
\begin{align}\label{r1}
\R{1}(x)=
\begin{cases}
t^*(x) \text{ if } KV_{0}(x) < \inf_{t \in \mathbb{R}^+} J(MV_{0},V_{0})(x,t), \\
\inf \{s \in \mathbb{R}^+ | J(MV_{0},V_{0})(x,s) \\
\quad < \inf_{t \in \mathbb{R}^+} J(MV_{0},V_{0})(x,t)+\epsilon\} \text{ otherwise,}
\end{cases}
\end{align}
and
\begin{align}\label{V1}
V_1(x)=
\begin{cases}
KV_{0}(x) \text{ if } KV_{0}(x) < \inf_{t \in \mathbb{R}^+} J(MV_{0},V_{0})(x,t), \\
J(MV_{0},V_{0})(x,\R{1}(x)) \text{ otherwise.}
\end{cases}
\end{align}
Then by induction, we define for $k\in\mathbb{N}^*$
\begin{align}\label{rk}
\R{k}(x)=
\begin{cases}
t^*(x) \text{ if } KV_{k-1}(x) < \inf_{t \in \mathbb{R}^+} J(MV_{k-1},V_{k-1})(x,t), \\
\inf \{s \in \mathbb{R}^+ | J(MV_{k-1},V_{k-1})(x,s) \\
< \inf_{t \in \mathbb{R}^+} J(MV_{k-1},V_{k-1})(x,t)+\epsilon\} \text{ otherwise}
\end{cases}
\end{align}
and
\begin{align}\label{vk}
V_k(x)=
\begin{cases}
KV_{k-1}(x) \text{ if } KV_{k-1}(x) \\
\quad < \inf_{t \in \mathbb{R}^+} J(MV_{k-1},V_{k-1})(x,t), \\
J(MV_{k-1},V_{k-1})(x,\R{k}(x)) \text{ otherwise.}
\end{cases}
\end{align}

\begin{rem}
As operators $K$ and $J(M\cdot,\cdot)$ send $\textbf{L}_{\Phi}(E)$ onto itself, see lemmas A.3 to A.6 and A.9 in \cite{BF2012}, and $h\in \textbf{L}_{\Phi}(E)$, one can prove that $V_k\in \textbf{L}_{\Phi}(E)$, for $k\in\mathbb{N}$.
\end{rem}


Now, we can show that the functions $V_k$ are close to the value function  $\mathcal{V}_k$ of our  optimal impulse control problem.

\begin{thm}\label{ineqVk}
Set $k\in\mathbb{N}$. Then we have the following inequality
\begin{equation}
\forall x \in E, \qquad \vc{k}(x) \leq V_{k}(x) \leq \vc{k}(x)+k \epsilon.
\end{equation}
\end{thm}

\textbf{Proof}  The case $k=1$ is straightforward because $\vc{0}=V_0=h$. 
For $k\geq 1$, we proceed by induction. We suppose that for some fixed $k$ one has
\[ \forall x \in E, \qquad \vc{k}(x) \leq V_{k}(x) \leq \vc{k}(x)+k \epsilon. \]
Applying lemma \ref{ineg_f&g} in the Appendix to this inequality, we obtain
\begin{equation}\label{inegK}
K\vc{k}(x) \leq KV_{k}(x) \leq K\vc{k}(x)+k\epsilon
\end{equation}
and for all $t \in \mathbb{R}^+$,
\begin{align}\label{inegJ}
J(M\vc{k},\vc{k})(x,t) \leq &J(MV_{k},V_{k})(x,t)  \nonumber \\
&\leq J(M\vc{k},\vc{k})(x,t)+k\epsilon.
\end{align}
We distinguish two cases. 

If $K\vc{k}(x) < \inf_{t \in \mathbb{R}^+} J(M\vc{k},\vc{k})(x,t)$, then by definition of $\vc{k+1}$, we have
\begin{equation}\label{ineq13}
 \vc{k+1}(x) =K\vc{k}(x)<\inf_{t \in \mathbb{R}^+} J(M\vc{k},\vc{k})(x,t).
 \end{equation}
The inequality \eqref{inegK} entails
\begin{equation}\label{inegK1}
\vc{k+1}(x) \leq KV_{k}(x) \leq \vc{k+1}(x)+k\epsilon.
\end{equation}
We have to distinguish the two cases arising in the definition of $V_{k+1}$.

$\bullet$ If $KV_{k}(x) < \inf_{t \in \mathbb{R}^+} J(MV_{k},V_{k})(x,t)$, then by definition \eqref{vk} of $V_{k+1}$, we have $V_{k+1}(x)=KV_{k}(x)$.
Therefore using inequality \eqref{inegK1}, we obtain the result:
\begin{equation*}
\vc{k+1}(x) \leq V_{k+1}(x) \leq \vc{k+1}(x)+ k\epsilon.
\end{equation*}
$\bullet$ We turn now to the case where
$\inf_{t \in \mathbb{R}^+} J(MV_{k},V_{k})(x,t)$ $\leq KV_{k}(x).$
By definition \eqref{vk} of $V_{k+1}$, we have
\begin{equation}\label{equ01}
V_{k+1}(x)=J(MV_{k},V_{k})(x,\R{k+1}(x)).
\end{equation}
On the one hand, inequality \eqref{inegJ} with $t=\R{k+1}(x)$ gives us
$J(M\vc{k},\vc{k})(x,\R{k+1}(x)) \leq J(MV_{k},V_{k})(x,\R{k+1}(x)).$ 
Thus we have
$J(M\vc{k},\vc{k})(x,\R{k+1}(x)) \leq V_{k+1}(x),$ using equation \eqref{equ01}.
We obtain
\[ \inf_{t \in \mathbb{R}^+} J(M\vc{k},\vc{k})(x,t) \leq V_{k+1}(x)\]
and the inequality \eqref{ineq13} yields
\begin{equation}\label{ineq05}
\vc{k+1}(x)\leq \inf_{t \in \mathbb{R}^+} J(M\vc{k},\vc{k})(x,t) \leq V_{k+1}(x).
\end{equation}
On the other hand, 
by definition \eqref{rk} of $\R{k+1}(x)$, and equality \eqref{equ01} we have
\begin{align*}
V_{k+1}(x)&=
J(MV_{k},V_{k})(x,\R{k+1}(x)) \\
&\leq  \inf_{t \in \mathbb{R}^+} J(MV_{k},V_{k})(x,t) + \epsilon \\
& \leq KV_{k}(x) + \epsilon,
\end{align*}
and the second part of inequality \eqref{inegK1} implies
\begin{equation}\label{ineq06}
V_{k+1}(x)  \leq \vc{k+1}(x) + (k+1)\epsilon.
\end{equation}
Combining inequalities \eqref{ineq05} and \eqref{ineq06}, we have the expected result
\begin{equation*}
\vc{k+1}(x)  \leq V_{k+1}(x)  \leq \vc{k+1}(x) + (k+1)\epsilon.
\end{equation*}
If $K\vc{k}(x) \geq \inf_{t \in \mathbb{R}^+} J(M\vc{k},\vc{k})(x,t)$, then by definition of $\vc{k+1}$ we have 
\begin{equation}\label{ineq07}
\vc{k+1}(x)=\inf_{t \in \mathbb{R}^+} J(M\vc{k},\vc{k})(x,t)\leq K\vc{k}(x).
\end{equation}
Again, we distinguish two sub-cases.

$\bullet$ If $KV_{k}(x) < \inf_{t \in \mathbb{R}^+} J(MV_{k},V_{k})(x,t)$, then by definition \eqref{vk} of $V_{k+1}$, we have
\begin{equation}\label{ineq09}
V_{k+1}(x)=KV_{k}(x)< \inf_{t \in \mathbb{R}^+} J(MV_{k},V_{k})(x,t) . 
\end{equation}
Next we use equations \eqref{inegJ} and \eqref{ineq07} to obtain
\begin{equation*}
V_{k+1}(x)  < \inf_{t \in \mathbb{R}^+} J(M\vc{k},\vc{k})(x,t)+ k\epsilon= \vc{k+1} (x)+k \epsilon . 
\end{equation*}
Using equations \eqref{ineq07}, \eqref{inegK} and \eqref{ineq09} then yields the expected result
\begin{eqnarray*}
\vc{k+1}(x) \leq K\vc{k}(x) \leq KV_{k}(x) &=& V_{k+1}(x)
\leq \vc{k+1}(x)+ k\epsilon.
\end{eqnarray*}
%
%
$\bullet$ If $KV_{k}(x) \geq \inf_{t \in \mathbb{R}^+} J(MV_{k},V_{k})(x,t)$, then  equation (\ref{equ01}) holds.
%
On the one hand, by definition \eqref{rk} of $\R{k+1}(x)$, we have
$ V_{k+1}(x)
\leq  \inf_{t \in \mathbb{R}^+} J(MV_{k},V_{k})(x,t) + \epsilon.$
Next we use equations \eqref{inegJ} and \eqref{ineq07} to obtain
\begin{align}\label{ineq11}
V_{k+1}(x)   \leq  \vc{k+1}(x)+ (k+1) \epsilon.
\end{align}
On the other hand, restarting with equality \eqref{ineq07} and using again inequality \eqref{inegJ}, we obtain
\begin{align*}
\vc{k+1}(x)
& \leq  \inf_{t \in \mathbb{R}^+} J(MV_{k},V_{k})(x,t) \\
&\leq  J(MV_{k},V_{k})(x,\R{k+1}(x)),
\end{align*}
and using equation \eqref{equ01} and inequality \eqref{ineq11}, we have
\[\vc{k+1}(x) \leq V_{k+1}(x) \leq \vc{k+1}(x)+ (k+1) \epsilon.\]
Finally, in all cases, we obtain the expected result.
\hfill $\Box$

\section{Construction of $\epsilon$-optimal strategies}\label{Sec4}

We now define an auxiliary PDMP on an enlarged state space. It can be interpreted as a controlled version of the original PDMP for a family of explicit strategies indexed by some integer $n$, where interventions are allowed up to the $n$-th jump time of the process. Then we establish that the cost of this strategy is exactly the approximate value function $V_n$.

\subsection{Construction of a controlled process $\{\T{X}_t\}$}\label{sect_def_strategie} 

Starting from the process $\{X_t\}$, we build a new PDMP $\{\T{X}_t\}$ on an enlarged state space. It has the same flow but a different jump mechanism from $\{X_t\}$.
The new process has four variables: $\T{X}_t=(m_t,\zeta_t,N_t,\theta_t)$, for $t\geq 0$. We consider that $\T{m}_t=(m_t,N_t)$ is the mode variable and $\T{\zeta}_t=(\zeta_t,\theta_t)$ is the Euclidean variable of the PDMP $\{\T{X}_t\}$. 
The state space of this process is defined by
\begin{align*}
\T{E}=\{(m,\zeta,N,\theta): (m,\zeta) \in E, N \in\mathbb{N}, \\
\theta \in[0;r_{\epsilon}^N(m,\zeta)\wedge t^*(m,\zeta)] \}\cup\{\Delta\},
\end{align*}
where $\Delta$ is an isolated point, called a cemetery state.
We define a distance on $\T{E}-\{\Delta\}$ in the same way as on $E$. For $(m,\zeta,N,\theta)$ and $(m',\zeta',N',\theta')$, setting $\T{m}=(m,N)$, $\T{m}'=(m',N')$, $\T{\zeta}=(\zeta,\theta)$ and $\T{\zeta}'=(\zeta',\theta')$, we define
\[\big|(\T{m},\T{\zeta})-(\T{m}',\T{\zeta}')\big|=\big|\T{\zeta}-\T{\zeta}' \big| \id{\T{m}=\T{m}'} + \infty \id{\T{m}\neq \T{m}'} . \] 
Let $(m,\zeta,N,\theta)\in\T{E}-\{\Delta\}$. Setting $\R{0}(m,\zeta)=+\infty$ and recalling that definition of $\R{N}$ is given by formulas \eqref{r1} and \eqref{rk}, we set $\T{r}_{\epsilon}(\Delta)=+\infty$ and
\begin{align}\label{defrepsilon}
\T{r}_{\epsilon}(m,\zeta,N,\theta)=
\R{N}(m,\zeta).
\end{align}
For $N\in \mathbb{N}^*$, define
\begin{eqnarray}\label{y_non_controlle}
\lefteqn{\y{N}(m,\zeta)}\nonumber\\
&=& \argmin_{y\in \mathbb{U}} \{ c(\Phi((m,\zeta),\R{N}(m,\zeta)),y)+V_{N-1}(y) \},
\end{eqnarray}
and set $\T{y}_\epsilon(m,\zeta,0,\theta)=  \Delta$ and if $N\in\mathbb{N}^*$
\begin{align}\label{y}
\T{y}_\epsilon(m,\zeta,N,\theta)= 
  (y^N_\epsilon(m,\zeta),N-1,0).
\end{align}
The local characteristics of the process are as follows. Consider $(m,\zeta,N,\theta)\in\T{E}-\{\Delta\}$.
\begin{itemize}
\item The flow $\T{\Phi}$ of process $\{\T{X}_t\}$ is  $\T{\Phi}(\Delta,t )=0$ and
\begin{equation}\label{phi_tilde}
\T{\Phi}((m,\zeta,N,\theta),t )=(m,\Phi_m(\zeta,t),N,\theta+t).
\end{equation}
The time the flow takes to reach the boundary of $\T{E}$ is denoted by $\T{t}^*$ and it is defined by
\begin{equation}\label{t_etoile_tilde}
\T{t}^*(m,\zeta,N,\theta)=t^*(m,\zeta)\wedge \R{N}(m,\zeta).
\end{equation}
\item The jump intensity $\T{\lambda}$ is defined by $\T{\lambda}(\Delta)=0$,
\begin{equation}\label{lambda_tilde}
\T{\lambda}(m,\zeta,N,\theta)=\lambda(m,\zeta),
\end{equation}
and setting
\[\T{\Lambda}(m,\zeta,N,\theta,t):=\int_0^t \T{\lambda}(m,\Phi_m(\zeta,s),N,\theta+s)ds,\]
we obtain
\begin{equation}\label{Lambda_tilde}
\T{\Lambda}(m,\zeta,N,\theta,t) = \Lambda(m,\zeta,t).
\end{equation}
\item Let $A\in \mathfrak{B}(\bar{E})$, $B\in\mathfrak{B}(\mathbb{N})$ and $C \in \mathfrak{B}(\mathbb{R}^+)$.
The Markov kernel $\T{Q}$ is defined by
\begin{align}\label{defqtilde}
\nonumber \T{Q}((m,\zeta&,N,\theta),A\times B\times C)\\
\nonumber &= Q((m,\zeta),A)\id{0}(C) \\
\nonumber &\times \Bigl[\Bigr.\id{N=0}\id{0}(B)+\id{N\geq 1}\id{N-1}(B) \\
\nonumber &\times \Bigl(\id{\theta<\R{N}(m,\zeta)} + \id{\theta=\R{N}(m,\zeta)=t^*(m,\zeta)}\Bigr)\Bigl.\Bigr] \\
 \nonumber&+\id{N\geq 1} \id{\theta=\R{N}(m,\zeta)\neq t^*(m,\zeta)}\\
 &\times\id{\T{y}_\epsilon(m,\zeta,N,\theta)} (A\times B\times C),
\end{align}
and $\T{Q}(\Delta,\Delta)=1$.
\end{itemize}
There exists a filtered probability space $(\T{\Omega}, \T{\mathcal{F}},(\T{\mathcal{F}})_{t\geq 0},$ $(\T{\p}_{\T{x}})_{\T{x} \in \T{E} } )$ on which the process $\{\T{X}_t\}$ is a strong Markov process \cite[sections 24-25]{Davis93}.

The additional components of the process $\{\T{X}_t\}$ can be interpreted as the following way: $\theta_t$ is the time since the last jump before $t$ and $N_t$ can be seen either as  a remaining time horizon or equivalently as a decreasing counter of jumps.
The process $\{\T{X}_t\}$ can be interpreted as a controlled version of $\{X_t\}$ where interventions are possible until the $N_0$-th jump. 
After the $N_0$-th jump, no more interventions are allowed and the process follows the same dynamics as $\{X_t\}$. Therefore, the control horizon $N_t$ decreases of one unit at each jump.

We denote by $\pi$ the projection of space $E\times \mathbb{N} \times \mathbb{R}$ on space $E$.
We denote by $(\T{S}_n)_{n\in\mathbb{N}^*}$ the sojourn times between two consecutive jumps, $(\T{Z}_n)_{n\in\mathbb{N}}$ the post jump locations and $(\T{T}_n)_{n\in\mathbb{N}}$ the jump times of the process $\{\T{X}_t\}$.
We distinguish two kinds of jump as follows, for $k\in \mathbb{N}^*$ :
\begin{itemize}
\item if $\T{S}_k < \T{r}_{\epsilon}(\T{Z}_{k-1})$ or if $\T{S}_k = \T{r}_{\epsilon}(\T{Z}_{k-1})=t^*(\pi(\T{Z}_{k-1}))$, then $\T{Z}_k$ is determined by the Markov kernel $Q$ of the non controlled process (equation \eqref{defqtilde}). We interpret the $k$-th jump as a \textit{natural jump}.
\item if $\T{S}_k = \T{r}_{\epsilon}(\T{Z}_{k-1}) \neq t^*(\pi(\T{Z}_{k-1}))$, then by equation \eqref{defqtilde}, we have $\T{Z}_k=\T{y}_{\epsilon}(\T{Z}_{k-1})$. We interpret the $k$-th jump as an \textit{intervention}. When an intervention occurs, the new starting point of $\{\T{X}_t\}$ is chosen by formulas \eqref{defqtilde} and \eqref{y}, where the set $\mathbb{U}$ is the control set. 
\end{itemize}

We now explicit the control strategy corresponding to the process $\{\T{X}_t\}$.
We introduce a counting process defined by
\[p^*(t):=\sum_{i=1}^\infty \id{\T{T}_i\leq t } \id{\T{X}_{\T{T}_i-} \in \partial \T{E}} \id{\T{S}_i\neq t^*(\pi(\T{Z}_{i-1}))}.\]
This process counts the number of jumps corresponding to interventions. Set
\begin{equation}\label{def_tps_int}
\T{\tau}_i= \inf \{t\in\mathbb{R}^+ \ | \ p^*(t)=i \}
\end{equation}
the time of the $i$-th intervention
and define the restarting points as
\begin{align*}
\T{R}_i= 
  \begin{cases}
  \T{Z}_{\T{\tau}_i} & \text{ if } \T{\tau}_i<\infty \\
  \Delta & \text{otherwise}.
  \end{cases}
\end{align*}
So the strategy associated to the controlled process $\{\T{X}_t\}$, starting at $(x,N_0,0)$ where $x\in E$ and $N_0\in\mathbb{N}$, is the sequence 
\[\T{\mathcal{S}}_\epsilon^{N_0}=(\T{\tau}_i,\T{R}_i)_{i\in\mathbb{N}^*}.\]
Ensuring that stategy $\T{\mathcal{S}}_\epsilon^{N_0}$ is admissible is a very difficult problem that will not be discussed here.



\subsection{Cost of the controlled trajectories}\label{sect_optimal}

We now turn to the calculation of the cost of the strategy $\T{\mathcal{S}}_\epsilon^{N_0}$ and compare it to the value function $\mathcal{V}_{N_0}$ 
More exactly, we  show that the cost of this strategy equals ${V}_{N_0}$.
To do so, we first need a technical result ensuring that $\R{k} (x) = t^*(x)$ if and only if $KV_{k-1}(x) < \inf_{t \in \mathbb{R}^+} J(MV_{k-1},V_{k-1})(x,t)$
\begin{prop} \label{R_epsilon}
Let $k\in\mathbb{N}^*$ and $x\in E$ such that $KV_{k-1}(x) \geq \inf_{t \in \mathbb{R}^+} J(MV_{k-1},V_{k-1})(x,t)$. Then we have
\[\R{k} (x) < t^*(x).\]
\end{prop} 

\textbf{Proof} 
By definition \eqref{def_J} of operator $J$ and definition \eqref{rk} of $\R{k}$, we have necessarily $\R{k}  (x) \leq t^*(x)$. 
Suppose $\R{k}  (x) =t^*(x)$. Then the definition of $\R{k}$ entails
\begin{align}\label{ineq-t-etoile}
\nonumber J(MV_{k-1},&V_{k-1})(x,t^*(x)) \\
& < \inf_{t \in \mathbb{R}^+} J(MV_{k-1},V_{k-1})(x,t)+\epsilon,
\end{align}
and for all $s \in [0;t^*(x)[,$
\begin{equation}\label{ineq-s}
J(MV_{k-1},V_{k-1})(x,s) \geq \inf_{t \in \mathbb{R}^+} J(MV_{k-1},V_{k-1})(x,t)+\epsilon.
\end{equation}
We define $d= \inf_{t \in \mathbb{R}^+} J(MV_{k-1},V_{k-1})(x,t)+\epsilon-J(MV_{k-1},V_{k-1})(x,t^*(x))$. The inequation \eqref{ineq-t-etoile} shows that $d>0$. By continuity of $J(MV_{k-1},V_{k-1})(x,.)$ on $\mathbb{R}^+$ (see proposition \ref{continuite_J} in the Appendix),
 there exists $\eta >0$ such that for all $s \in [t^*(x)-\eta;t^*(x)[$
\begin{align}\label{ineq-abs}
|J(MV_{k-1},V_{k-1})(x,s)
-J(MV_{k-1},V_{k-1})(x,t^*(x))|<d.
\end{align}
Consider $s\in [t^*(x)-\eta;t^*(x)[$. Two cases are possible.\\
$\bullet$ If $J(MV_{k-1},V_{k-1})(x,s)\leq J(MV_{k-1},V_{k-1})(x,t^*(x))$, then by inequation \eqref{ineq-t-etoile} we have
\[J(MV_{k-1},V_{k-1})(x,s) < \inf_{t \in \mathbb{R}^+} J(MV_{k-1},V_{k-1})(x,t)+\epsilon \]
contradicting the inequation \eqref{ineq-s}.\\
$\bullet$ If $J(MV_{k-1},V_{k-1})(x,s) > J(MV_{k-1},V_{k-1})(x,t^*(x))$, then by \eqref{ineq-abs} and definition of $d$, we obtain
\begin{align*}
J(MV_{k-1},V_{k-1})(x,s) < &\inf_{t \in \mathbb{R}^+} J(MV_{k-1},V_{k-1})(x,t)+\epsilon
\end{align*}
contradicting again the inequation \eqref{ineq-s}.
We conclude that $\R{k} (x) \neq t^*(x)$ in both cases.
\hfill $\Box$


We define the following cost functions of strategies $\tilde{\mathcal{S}}^{N_0}_\epsilon$, $N_0\in\mathbb{N}$. If $N_0=0$,
\begin{equation}\label{cout_J0}
\mathcal{J}_0(x)=\et{0}\Bigl[ \int_0^{\infty} \e^{-\alpha s}\T{f}(\T{X}_s)ds \Bigr],
\end{equation}
and for $N_0 \geq 1$,
\begin{align}\label{cout_JN0}
\nonumber \mathcal{J}_{N_0}(x)=\et{N_0}\Bigl[\Bigr. \int_0^{\infty} \e^{-\alpha s}\T{f}(\T{X}_s)ds \\
+ \sum_{i=1}^{\infty} \id{\T{\tau}_i<\infty}\e^{-\alpha \T{\tau}_i } \T{c}(\T{X}_{\T{\tau}_i-},\T{X}_{\T{\tau}_i})\Bigl.\Bigr],
\end{align}
where 
functions $\T{f}$ and $\T{c}$ are defined by
\begin{align}\label{f_c_tilde}
\forall t \in \mathbb{R}^+, \quad 
\begin{cases}\T{c}(\T{X}_{t-},\T{X}_t)=c(\pi(\T{X}_{t-}),\pi(\T{X}_t)) \\
\T{f}(\T{X}_t)=f(\pi(\T{X}_t)),
\end{cases}
\end{align}
recalling that $\pi$ is the projection of space $\T{E}$ on space $E$.
Note that the sum in the expectation is always finite.

We now prove by induction that $\mathcal{J}_n=V_n$ for all $n\geq 0$. The proof is split in the three following theorems. We first prove that $\mathcal{J}_0=h=V_0$ and then that the sequence $(\mathcal{J}_n)$ satisfies the same recursion as the sequence $(V_n)$.

\begin{thm}\label{cost_strategy_0}
Set $x\in E$. We have the following equality:
\begin{equation}\label{j0}
\mathcal{J}_0(x)=h(x)=\E_x \Bigl[ \int_0^{\infty} \e^{-\alpha s}f(X_s)ds \Bigr].
\end{equation}
\end{thm}

\textbf{Proof}  Set $x=(m,\zeta)\in E$.
By definition of $\mathcal{J}_0$ and $\T{f}$, we have
\begin{align*}
\mathcal{J}_0(x)
=\et{0}\Bigl[ \int_0^{\infty} \e^{-\alpha s}f(\pi(\T{X}_s))ds \Bigr].
\end{align*}
We show that $\{\pi(\T{X}_t)\}$ starting at $\T{x}=(x,0,0)$ is a PDMP which has the same characteristics as the process $\{X_t\}$ starting at $x$. \\
The process $\{\pi(\T{X}_t)\}$ starts at the point $\pi(\T{X}_0)=(m,\zeta)=x$.
It then follows the flow (see \eqref{phi_tilde}) 
\[ \pi(\T{\Phi}((m,\zeta,0,\theta),t))=(m,\Phi_m(\zeta,t))=\Phi(x,t) \]
until the first jump time, which is determined by the survival equation:
\[ \T{\p}_{\T{x}}(\T{T}_1 > t)= \exp(-\T{\Lambda}(\T{x})) \id{t<\T{t}^*(\T{x})}.\]
By definition \eqref{t_etoile_tilde} of $\T{t}^*$ and \eqref{defrepsilon} of $\T{r}_\epsilon$, we have 
$\T{t}^*(\T{x})=t^*(x)$. 
Furthermore, expression \eqref{Lambda_tilde} of $\T{\Lambda}$ gives us
\begin{align*}
\T{\p}_{\T{x}}(\T{T}_1 > t)= \exp(-\Lambda(x)) \id{t<t^*(x)} 
=\p_{x}(T_1 > t).
\end{align*}
At this random time, a jump occurs and the process $\{\pi(\T{X}_t)\}$ restarts at a new point $\pi(\T{X}_{\T{T}_1})=(m_{\T{T}_1},\zeta_{\T{T}_1})$ according to the distribution (see \eqref{defqtilde})
\begin{align*}
\T{Q}(\T{\Phi}(\T{x},\T{T}_1),(\cdot,0,0))
=&\T{Q}((\Phi(x,\T{T}_1),0,\theta_{\T{T}_1}),(\cdot,0,0)) \\
=& Q(\Phi(x,\T{T}_1),\cdot) \\
=& Q(\Phi(x,T_1),\cdot).
\end{align*}
Then we can continue the construction of the process $\{\pi(\T{X}_t)\}$ in a similar way. We conclude that $\{\pi(\T{X}_t)\}$ is a PDMP with the same characteristics as the PDMP $\{X_t\}$.
The consequence is that $\mathcal{J}_0$ is equal to $h$ (see definition \eqref{def_h}).
\hfill $\Box$

\begin{thm}\label{cost_strategy_1}
Set $x\in E$. We have the following equality:
\begin{align*}
\mathcal{J}_1(x)=
  \begin{cases}
    K\mathcal{J}_0(x) \text{ if } K\mathcal{J}_0(x) < \inf_{t \in \mathbb{R}^+} J(M\mathcal{J}_0,\mathcal{J}_0)(x,t), \\
    J(M\mathcal{J}_0,\mathcal{J}_0)(x,\R{1}(x))\text{ otherwise.}
  \end{cases}
  \end{align*}
implying $\mathcal{J}_1=V_1$ because $\mathcal{J}_0=h$.
\end{thm}

\textbf{Proof}  Set $x\in E$ and $\T{x}:=(x,1,0)$ and recall that $h=Kh$ by \cite[Proposition 1]{Costa89}. There are again two cases.

If $K\mathcal{J}_0(x) < \inf_{t \in \mathbb{R}^+} J(M\mathcal{J}_0,\mathcal{J}_0)(x,t)$, 
then by theorem \ref{cost_strategy_0} we have $h(x)=Kh(x) < \inf_{t \in \mathbb{R}^+} J(Mh,h)(x,t)$. Thus formula \eqref{r1} gives $\R1(x)=t^*(x)$, so 
that the next jump is a natural jump. We are going to prove that $\{\pi(\T{X}_t)\}$ starting at $\T{x}=(x,1,0)$ has the same distribution as $\{X_t\}$ starting at $x$. \\
Formula \eqref{cout_JN0} for the cost, formulas \eqref{defqtilde} of $\T{Q}$, \eqref{f_c_tilde} of $\T{f}$ and \eqref{phi_tilde} of $\T{\Phi}$ give
\begin{align*}
\mathcal{J}_1(x)
=&\et{1} \Bigl[ \int_0^{\T{S}_1} \e^{-\alpha s}f(\Phi(x,s))ds \Bigr. \\
& + \e^{-\alpha \T{S}_1} \et{1} \Bigl[\int_{0}^{\infty} \e^{-\alpha t}\T{f}(\T{X}_{t+\T{S}_1})dt \Big| \T{\mathcal{F}}_{\T{S}_1}  \Bigr]\Bigl.\Bigr].
 \end{align*}
The Markov property for $\{\T{X}_t\}$ entails
\begin{align*}
\mathcal{J}_1(x)
=&\et{1} \Bigl[ \int_0^{\T{S}_1} \e^{-\alpha s}f(\Phi(x,s))ds \Bigr] \\
&+ \et{1} \Bigl[  \e^{-\alpha \T{S}_1} 
\T{E}_{\T{Z}_1}\Bigl[
\int_{0}^{\infty} \e^{-\alpha t}\T{f}(\T{X}_t)dt 
\Bigl]
\Bigr].
 \end{align*}
As $N_0=1$ and $t^*(x)=\R{1}(x)$, we have either $S_1<t^*(x)\wedge\R{1}(x)$ or $S_1=t^*(x)=\R{1}(x)$. Then by Lemma \ref{loi_conjointe} in the Appendix, we obtain
\begin{align*}
\mathcal{J}_1(x)
=&\E_x \Bigl[ \int_0^{S_1} \e^{-\alpha s}f(\Phi(x,s))ds \Bigr] \\
&+ \E_x \Bigl[  \e^{-\alpha S_1} 
\T{\E}_{(Z_1,0,0)} \Bigl[ \int_{0}^{\infty} \e^{-\alpha t}\T{f}(\T{X}_t)dt \Bigr] \Bigr]. 
\end{align*}
Definition \eqref{cout_J0} of $\mathcal{J}_0$ and definition \eqref{def_F} of $F$ imply
\begin{align*}
\mathcal{J}_1(x)
=&F(x,t^*(x))
+ \E_x \Bigl[  \e^{-\alpha S_1} \mathcal{J}_0(Z_1) \Bigr]. 
 \end{align*}
Finally definition \eqref{def_K} of operator $K$ and theorem \ref{cost_strategy_0} give us
$\mathcal{J}_1(x)= K \mathcal{J}_0(x)=Kh(x)$.

If $K\mathcal{J}_0(x) \geq \inf_{t \in \mathbb{R}^+} J(M\mathcal{J}_0,\mathcal{J}_0)(x,t)$, then by theorem \ref{cost_strategy_0}, we have 
$Kh(x) \geq \inf_{t \in \mathbb{R}^+} J(Mh,h)(x,t)$. Thus proposition \ref{R_epsilon} entails that $\R1(x) < t^*(x)$. Particularly, we have
$\R1(x) \neq t^*(x)$. The expression \eqref{defqtilde} of $\T{Q}$ shows that an intervention is possible. \\
First, we can notice that $\T{S}_1=\R{1}(x)$ is equivalent to $\T{S}_1\geq \R{1}(x)$ because we have always $\T{S_1} \leq \T{t}^*(\T{x})=\R{1}(x)$ by definition \eqref{t_etoile_tilde} of $\T{t}^*(\T{x})$. Then, starting from the definition \eqref{cout_JN0} of $\mathcal{J}_{N_0}$ for $N_0=1$ and using definitions \eqref{f_c_tilde} of $\T{c}$ and \eqref{phi_tilde} of $\T{\Phi}$, we have
\begin{equation}\label{decomposition}
\mathcal{J}_1(x) =\phi_1(x)+\phi_2(x)+\phi_3(x).
\end{equation}
where
\begin{align*}
\phi_1(x):= &\et{1} \Bigl[ \int_{0}^{\T{S}_1\wedge \R1(x)}\e^{-\alpha s}\T{f}(\T{X}_s)ds \Bigr] \\
\phi_2(x):= &\et{1} \Bigl[ \id{\T{S}_1< \R1(x)} \int_{\T{S}_1}^{+\infty} \e^{-\alpha s}\T{f}(\T{X}_s)ds \Bigr] \\
\phi_3(x):= &\et{1} \Bigl[\Bigr. \id{\T{S}_1\geq \R1(x)} \Big\{ \int_{\R1(x)}^{+\infty} \e^{-\alpha s}\T{f}(\T{X}_s)ds \\
& + \e^{-\alpha \R1(x)}c(\Phi(x,\R1(x)),\y1(x)) \Big\} \Bigl.\Bigr].
\end{align*}
\textbf{For the first part of \eqref{decomposition}}, we have
\begin{align*}
\phi_1(x)
= & \et{1} \Bigl[ \int_{0}^{\T{S}_1\wedge \R1(x)}\e^{-\alpha s} f(\Phi(x,s))ds \Bigr].
\end{align*}
Recalling that $N_0=1$ and $\R{1}(x)\neq t^*(x)$, we apply lemma \ref{loi_conjointe} and we obtain
\begin{align*}
\phi_1(x)
=& \E_x \Bigl[ \int_{0}^{S_1\wedge \R1(x)}\e^{-\alpha s} f(\Phi(x,s))ds \Bigr].
\end{align*}
Finally, definition \eqref{def_F} of $F$ entails
\begin{equation}\label{phi1}
\phi_1(x) = F(x,r^1_\epsilon(x)).
\end{equation}
\textbf{For the second part of \eqref{decomposition}}, we have
\begin{align*}
\phi_2(x)
= &\et{1} \Bigl[\Bigr. \id{\T{S}_1< \R1(x)} \e^{-\alpha \T{S}_1} \\
& \et{1} \Bigl[\int_{0}^{+\infty} \e^{-\alpha t}\T{f}(\T{X}_{t+\T{S}_1})dt \Big| \T{\mathcal{F}}_{\T{S}_1} \Bigr] \Bigl.\Bigr].
\end{align*}
The Markov property implies
\begin{align*}
\phi_2(x)
= &\et{1} \Bigl[\Bigr. \id{\T{S}_1< \R1(x)} \e^{-\alpha \T{S}_1} \\
&\times
  \T{\E}_{\T{Z}_1} \Bigl[\int_{0}^{+\infty} \e^{-\alpha t}\T{f}(\T{X}_{t})dt \Bigr]\Bigl.\Bigr].
\end{align*}
Recalling that $N_0=1$ and $\R{1}(x)\neq t^*(x)$, we apply again lemma \ref{loi_conjointe}
\begin{align*}
\phi_2(x)
= &\E_x \Bigl[\Bigr. \id{S_1< \R1(x)} \e^{-\alpha S_1} \\
&\times  \T{\E}_{(Z_1,0,0)} \Bigl[\int_{0}^{+\infty} \e^{-\alpha t}\T{f}(\T{X}_{t})dt \Bigr]\Bigl.\Bigr].
\end{align*}
Definition \eqref{cout_J0} of $\mathcal{J}_0$  and after theorem \ref{cost_strategy_0} entail
\begin{equation}\label{phi2}
\phi_2(x) = \E_x \left[  \id{S_1< \R1(x)}  \e^{-\alpha S_1} h(Z_1)\right].
\end{equation}
\textbf{For the third part of \eqref{decomposition}}, we have
\begin{align*}
\phi_3(x)
 = &\et{1} \Bigl[ \Bigr. \id{\T{S}_1\geq \R1(x)}  \e^{-\alpha \R1(x)} \Big\{ \\
 &\et{1}\Bigl[\int_{0}^{+\infty} \e^{-\alpha t}\T{f}(\T{X}_{t+\R1(x)})dt \Big| \T{\mathcal{F}}_{\R1(x)}\Bigr] \\
  & +c(\Phi(x,\R1(x)),\y1(x)) \Big\} \Bigl. \Bigr]. 
\end{align*}
By the Markov property, we obtain
\begin{align*}
\phi_3(x)
 = &\et{1} \Bigl[ \Bigr. \id{\T{S}_1\geq \R1(x)}  \e^{-\alpha \R1(x)} \Big\{ \\
& \T{\E}_{(y^1_\epsilon(x),0,0)} \Bigl[\int_{0}^{+\infty} \e^{-\alpha t}\T{f}(\T{X}_{t})dt \Bigr]  \\
& +c(\Phi(x,\R1(x)),\y1(x)) \Big\} \Bigl. \Bigr]. 
\end{align*}
First using definition \eqref{cout_J0} of $\mathcal{J}_0$ and theorem \ref{cost_strategy_0}, and after using definition \eqref{def_M} of operator $M$ and definition \eqref{y} of $y_\epsilon^1$, we have
\begin{align*}
\phi_3(x)
 = &\et{1} \Bigl[ \Bigr. \id{\T{S}_1\geq \R1(x)}  \e^{-\alpha \R1(x)} \Big\{ h(y^1_\epsilon(x))  \\
 & +c(\Phi(x,\R1(x)),\y1(x)) \Big\} \Bigl. \Bigr] \\
 = & \et{1} \Bigl[ \Bigr. \id{\T{S}_1\geq \R1(x)}  \e^{-\alpha \R1(x)} Mh(\Phi(x,\R1(x))) \Bigl. \Bigr].
\end{align*}
Recalling that $N_0=1$ and $\R{1}(x)\neq t^*(x)$, we apply again lemma \ref{loi_conjointe} and we obtain
\begin{equation}\label{phi3}
\phi_3(x) = \E_x\Bigl[  \e^{-\alpha \R1(x)} Mh(\Phi(x,\R1(x))) \id{S_1\geq \R1(x)} \Bigr].
\end{equation}
Finally, the equations \eqref{decomposition}, \eqref{phi1}, \eqref{phi2} and \eqref{phi3} yield
\begin{align*}
\mathcal{J}_1(x)=& F(x,\R1(x)) +  \E_x \Bigl[\Bigr. \e^{-\alpha \R1(x)} Mh(\Phi(x,\R1(x))) \\
&\id{S_1\geq \R1(x)} + \e^{-\alpha S_1} h(Z_1)  \id{S_1< \R1(x)} \Bigl.\Bigr]
\end{align*}
and definition \eqref{def_J} of operator $J$ gives us the result $\mathcal{J}_1(x)= J(Mh,h)(x,\R1(x))$.
\hfill $\Box$

\begin{thm}\label{cost_strategy_n}
Let $N_0$ be an integer such that $N_0 \geq 1$. Let $x\in E$. We have the formula
\begin{align*}
\mathcal{J}_{N_0}(x)= 
  \begin{cases}
    K\mathcal{J}_{N_0-1}(x) \text{ if } K\mathcal{J}_{N_0-1}(x)< \\
     \qquad \qquad \inf_{t \in \mathbb{R}^+} J(M\mathcal{J}_{N_0-1},\mathcal{J}_{N_0-1})(x,t), \\
    J(M\mathcal{J}_{N_0-1},\mathcal{J}_{N_0-1})(x,\R{N_0}(x)) \text{ otherwise.}
  \end{cases}
  \end{align*}
  Furthermore, we have $\mathcal{J}_{N_0}=V_{N_0}$.
\end{thm}


\textbf{Proof}   Some details are similar to the previous proof and are therefore omitted. We proceed by induction. The case $N_0=1$ was showed before in theorem \ref{cost_strategy_1}. Let $N_0 \geq 1$ and $x\in E$. Suppose that 
\begin{align*}
\mathcal{J}_{N_0}(x)= 
  \begin{cases}
    K\mathcal{J}_{N_0-1}(x) \text{ if } K\mathcal{J}_{N_0-1}(x) \\
    \qquad < \inf_{t \in \mathbb{R}^+} J(M\mathcal{J}_{N_0-1},\mathcal{J}_{N_0-1})(x,t), \\
    J(M\mathcal{J}_{N_0-1},\mathcal{J}_{N_0-1})(x,\R{N_0}(x)) \quad  \text{otherwise.}
  \end{cases}
  \end{align*}
and $\mathcal{J}_{N_0}=V_{N_0}$. 
%
We again distinguish two cases.

If $K\mathcal{J}_{N_0}(x)< \inf_{t \in \mathbb{R}^+} J(M\mathcal{J}_{N_0},\mathcal{J}_{N_0})(x,t)$, then by the induction hypothesis, we obtain $KV_{N_0}(x)< \inf_{t \in \mathbb{R}^+} J(MV_{N_0},V_{N_0})(x,t)$. Whence by definition \eqref{rk} of $\R{N_0+1}$, we have 
$\R{N_0+1}(x)=t^*(x)$. 
Thus definition \eqref{defqtilde} of $\T{Q}$ show that the first jump is not an intervention. Consequently, we have
\begin{align*}
\mathcal{J}_{N_0+1}(x)
=& \psi_1(x)+\psi_2(x), 
\end{align*}
where
\begin{align*}
\psi_1(x):=&\et{N_0+1} \Bigl[ \int_0^{\T{S}_1} \e^{-\alpha s} \T{f}(\T{\Phi}((x,N_0+1,0),s)) ds \Bigr] \\
\psi_2(x):=&\et{N_0+1} \Bigl[\Bigr. \int_{\T{S}_1}^\infty \e^{-\alpha s} \T{f}(\T{X}_s) ds \\
&+\sum_{i=1}^{\infty} \e^{-\alpha \T{\tau}_i } \T{c}(\T{X}_{\T{\tau}_i -},\T{X}_{\T{\tau}_i}) \Bigl.\Bigr] . 
\end{align*}
\textbf{For the first part of $\mathcal{J}_{N_0+1}$}, using lemma \ref{loi_conjointe} with $N_0\geq 1$ and $\R{N_0+1}(x)=t^*(x)$, we have
\begin{equation}\label{psi1}
\psi_1(x)=F(x,t^*(x)).
\end{equation}
\textbf{For the second part of $\mathcal{J}_{N_0+1}$}, we have
\begin{align*}
\psi_2(x)
 =&\et{N_0+1} \Bigl[ \e^{-\alpha \T{S}_1} 
  \et{N_0+1} \Bigl[\Bigr.  \int_{0}^\infty \e^{-\alpha t} \T{f}(\T{X}_{t+\T{S}_1}) dt \\
 &+ \sum_{i=1}^{\infty} \e^{-\alpha (\T{\tau}_i -\T{S}_1)} \T{c}(\T{X}_{\T{\tau}_i -},\T{X}_{\T{\tau}_i }) \Big| \T{\mathcal{F}}_{\T{S}_1} \Bigl.\Bigr] \Bigr] .
\end{align*}
Recalling that $N_0\geq 1$ and $\R{N_0+1}(x)=t^*(x)$, by proposition \ref{prop_markov_tau}, Markov property for process $\{\T{X}_t\}$ and lemma \ref{loi_conjointe}, we obtain
\begin{align*}
\psi_2(x)
=&\E_x \Bigl[\Bigr. \e^{-\alpha S_1} 
  \T{\E}_{(Z_1,N_0,0)} \Bigl[ \Bigr. \int_{0}^\infty \e^{-\alpha t} \T{f}(\T{X}_t) dt \\
 & + \sum_{i=1}^{\infty} \e^{-\alpha \T{\tau}_i } \T{c}(\T{X}_{\T{\tau}_i -},\T{X}_{\T{\tau}_i }) \Bigl.\Bigr] \Bigl.\Bigr]. 
\end{align*}
By definition \eqref{cout_JN0} of $\mathcal{J}_{N_0}$, we finally have
\begin{equation}\label{psi2}
\psi_2(x)=\E_x \Bigr[ \e^{-\alpha S_1} \mathcal{J}_{N_0}(Z_1) \Bigl].
\end{equation}
Thus equations \eqref{psi1} and \eqref{psi2} and definition \eqref{def_K} of operator $K$ give us $\mathcal{J}_{N_0+1}(x)= K\mathcal{J}_{N_0}(x)$.

If $K\mathcal{J}_{N_0}(x) \geq \inf_{t \in \mathbb{R}^+} J(M\mathcal{J}_{N_0},\mathcal{J}_{N_0})(x,t)$, then by the induction hypothesis, we have $KV_{N_0}(x)\geq \inf_{t \in \mathbb{R}^+} J(MV_{N_0},V_{N_0})(x,t)$. Consequently, proposition \ref{R_epsilon} entails $\R{N_0+1}(x)< t^*(x)$ and we have
$\R{N_0+1}(x)\neq t^*(x)$. \\
To calculate the cost of strategy $\mathcal{S}_\epsilon^{N_0+1}$ in this case, we start from formula \eqref{cout_JN0}.
Definition \eqref{defqtilde} of $\T{Q}$ shows that the first jump of the controlled process $\{\T{X}_t\}$ can be an intervention. It gives us
\begin{equation}\label{sum_Psi}
\mathcal{J}_{N_0+1}(x)=\Psi_1(x)+\Psi_2(x)+\Psi_3(x)+\Psi_4(x)
\end{equation}
where
\begin{align*}
\Psi_1(x):= &\et{N_0+1}\Bigl[ \int_0^{\T{S}_1\wedge \R{N_0+1}(x)} \e^{-\alpha s}\T{f}(\T{X}_s)ds\Bigr]  \\
\Psi_2(x):=& \et{N_0+1}\Bigl[\Bigr. \id{\T{S}_1< \R{N_0+1}(x)} \Big\{\int_{\T{S}_1}^\infty \e^{-\alpha s}\T{f}(\T{X}_s)ds \\
&+ \sum_{i=1}^{+\infty} \e^{-\alpha \T{\tau}_i } \T{c}(\T{X}_{\T{\tau}_i-},\T{X}_{\T{\tau}_i})\Big\} \Bigl.\Bigr] \\
\Psi_3(x):=& \et{N_0+1}\Bigl[\Bigr. \id{\T{S}_1= \R{N_0+1}(x)}  \e^{-\alpha \R{N_0+1}(x)} \\
&\T{c}( \T{\Phi}(\T{x},\R{N_0+1}(x)),(y^{N_0+1}_\epsilon(x),N_0,0))  \Bigl.\Bigr] \\
\Psi_4(x):=&  \et{N_0+1}\Bigl[\Bigr. \id{\T{S}_1= \R{N_0+1}(x)} \Big\{ \int_{\R{N_0+1}(x)}^\infty \e^{-\alpha s} \\
&\T{f}(\T{X}_s)ds + \sum_{i=2}^{+\infty} \e^{-\alpha \T{\tau}_i } \T{c}(\T{X}_{\T{\tau}_i-},\T{X}_{\T{\tau}_i}) \Big\} \Bigl.\Bigr].
\end{align*}
\textbf{For the first part of $\mathcal{J}_{N_0+1}$}, using lemma \ref{loi_conjointe} with $N_0\geq1$ and $\R{N_0+1}(x)<t^*(x)$, we obtain
\begin{equation}\label{Psi1}
\Psi_1(x)=F(x,\R{N_0+1}(x)).
\end{equation}
\textbf{For the second part of $\mathcal{J}_{N_0+1}$}, we have:
\begin{align*}
\Psi_2(x)
 =& \et{N_0+1}\Bigl[\Bigr. \id{\T{S}_1< \R{N_0+1}(x)} \e^{-\alpha\T{S}_1} \\
  & \et{N_0+1}\Bigl[\Bigr. \int_{0}^\infty \e^{-\alpha t}\T{f}(\T{X}_{t+\T{S}_1})dt  \\
 &+ \sum_{i=1}^{+\infty} \e^{-\alpha( \T{\tau}_i -\T{S}_1)} \T{c}(\T{X}_{\T{\tau}_i-},\T{X}_{\T{\tau}_i})\Big| \T{\mathcal{F}}_{\T{S}_1} \Bigl.\Bigr] \Bigl.\Bigr]. 
\end{align*}
Recalling that $N_0\geq 1$ and $\R{N_0+1}\neq t^*(x)$, by proposition \ref{prop_markov_tau} (with $\R{N_0+1}(x)<t^*(x)$), Markov property for process $\{\T{X}_t\}$ and lemma \ref{loi_conjointe}, we obtain
\begin{align*}
\Psi_2(x)
 =& \E_x\Bigl[\Bigr. \id{S_1< \R{N_0+1}(x)} \e^{-\alpha S_1} \\
  & \T{\E}_{(Z_1,N_0,0)} \Bigl[\Bigr. \int_{0}^\infty \e^{-\alpha t}\T{f}(\T{X}_{t})dt  \\
 &+ \sum_{i=1}^{+\infty} \e^{-\alpha \T{\tau}_i } \T{c}(\T{X}_{\T{\tau}_i-},\T{X}_{\T{\tau}_i}) \Bigl.\Bigr] \Bigl.\Bigr]. 
\end{align*}
Finally, definition \eqref{cout_JN0} of $\mathcal{J}_{N_0}$ gives us
\begin{equation}\label{Psi2}
\Psi_2(x) = \E_x \Bigl[ \id{S_1< \R{N_0+1}(x)} \e^{-\alpha S_1}  \mathcal{J}_{N_0}(Z_1) \Bigr].
\end{equation}
\textbf{For the third part of $\mathcal{J}_{N_0+1}$}, we can remark $\{\T{S}_1= \R{N_0+1}(x)\}$ is the same event as $\{\T{S}_1\geq \R{N_0+1}(x)\}$ because we have always $\T{S}_1 \leq \T{t}^*(\T{x})=\R{N_0+1}(x)$. 
Using lemma \ref{loi_conjointe} with $N_0\geq 1$ and $\R{N_0+1}\neq t^*(x)$, we can conclude
\begin{align}\label{Psi3}
\nonumber \Psi_3(x) &= \E_x \Bigl[\Bigr. \id{S_1\geq \R{N_0+1}(x)}\e^{-\alpha \R{N_0+1}(x)} \\
&c( \Phi(x,\R{N_0+1}(x)),y^{N_0+1}_\epsilon(x)) \Bigl.\Bigr]. 
\end{align}
\textbf{For the fourth part of $\mathcal{J}_{N_0+1}$}, we have
\begin{align*}
\Psi_4(x)
=& \et{N_0+1}\Bigl[\Bigr. \id{\T{S}_1= \R{N_0+1}(x)} \e^{-\alpha \R{N_0+1}(x)} \\
& \et{N_0+1}\Bigl[\Bigr.  \int_{0}^\infty \e^{-\alpha t}\T{f}(\T{X}_{t+\R{N_0+1}(x)})dt \\
&+ \sum_{i=2}^{+\infty} \e^{-\alpha (\T{\tau}_i -\R{N_0+1}(x))} \T{c}(\T{X}_{\T{\tau}_i-},\T{X}_{\T{\tau}_i}) 
\Big| \T{\mathcal{F}}_{\R{N_0+1}(x)} \Bigl.\Bigr] \Bigl.\Bigr] .
\end{align*}
Recalling $N_0\geq1$ and $\R{N_0+1}\neq t^*(x)$, by proposition \ref{prop_markov_tau} (with $\R{N_0+1}(x)<t^*(x)$), Markov property for process $\{\T{X}_t\}$ and lemma \ref{loi_conjointe}, we have
\begin{align*}
\Psi_4(x)
=& \E_x \Bigl[\Bigr. \id{S_1= \R{N_0+1}(x)} \e^{-\alpha \R{N_0+1}(x)} \T{\E}_{(Z_1,N_0,0)} \Bigl[\Bigr. \\
& \int_{0}^\infty \e^{-\alpha t}\T{f}(\T{X}_t)dt + \sum_{i=1}^{+\infty} \e^{-\alpha \T{\tau}_i } \T{c}(\T{X}_{\T{\tau}_i-},\T{X}_{\T{\tau}_i}) 
\Bigl.\Bigr] \Bigl.\Bigr]. 
\end{align*}
Using definition \eqref{cout_JN0} of $\mathcal{J}_{N_0}$
and by the same argument as the third part of $\mathcal{J}_{N_0+1}$, we finally obtain
\begin{equation}\label{Psi4}
\Psi_4(x)= \E_x \Bigl[ \id{S_1 \geq \R{N_0+1}(x)} \e^{-\alpha \R{N_0+1}(x)} \mathcal{J}_{N_0}(y_\epsilon^{N_0+1}(x)) \Bigr].
\end{equation}

Now, if we combine equations \eqref{sum_Psi} \eqref{Psi1}, \eqref{Psi2}, \eqref{Psi3} and \eqref{Psi4}, the expression of the cost $\mathcal{J}_{N_0+1}$ becomes
\begin{align*}
\mathcal{J}_{N_0+1}(x) 
&=F(x,\R{N_0+1}(x)) + \E_x \Bigl[\Bigr. \id{S_1< \R{N_0+1}(x)} \e^{-\alpha S_1}  \\
&\mathcal{J}_{N_0}(Z_1) + \id{S_1\geq \R{N_0+1}(x)}\e^{-\alpha \R{N_0+1}(x)}\\
& \Big\{ c( \Phi(x,\R{N_0+1}(x)),y^{N_0+1}_\epsilon(x)) + \mathcal{J}_{N_0}(y_\epsilon^{N_0+1}(x)) \Big\} \Bigl.\Bigr].
\end{align*}
Definition \eqref{def_M} of operator $M$ and definition \eqref{def_J} of operator $J$ entail
\begin{align*}
\mathcal{J}_{N_0+1}(x)= J(M\mathcal{J}_{N_0},\mathcal{J}_{N_0})(x,\R{N_0+1}(x)).
\end{align*}

By induction hypothesis, we have: $\mathcal{J}_{N_0}=V_{N_0}$. Then $\mathcal{J}_{N_0+1}$ can be written by
\begin{align*}
\mathcal{J}_{N_0+1}(x)=
  \begin{cases}
    KV_{N_0}(x) \text{ if } KV_{N_0}(x) < \\
    \qquad \qquad  \inf_{t \in \mathbb{R}^+} J(MV_{N_0},V_{N_0})(x,t), \\
    J(MV_{N_0},V_{N_0})(x,\R{N_0+1}(x))\text{ otherwise.}
  \end{cases}
  \end{align*}
corresponding exactly to the definition \eqref{vk} of $V_{N_0+1}$. Consequently we have $\mathcal{J}_{N_0+1}=V_{N_0+1}$.
\hfill $\Box$

%

From theorems \ref{ineqVk} and \ref{cost_strategy_n}, we thus infer that the family of strategies $\T{\mathcal{S}}_\epsilon^{N_0}$ have cost functions that are arbitrarily close to the value function when the parameter $\epsilon$ goes to $0$. Hence, we have explicitly constructed a family of $\epsilon$-optimal strategies for the impulse control problem under study. In a further work, we intend to propose a computable numerical approximation of such strategies, based on the ideas of \cite{BF2012}.

\appendix

\section{Technical results about functions $V_k$}

We first prove a result concerning the continuity of function $J(MV_{k},V_k)(x,\cdot)$ on $\mathbb{R}^+$. 

\begin{prop}\label{continuite_J}
For all $k\in\mathbb{N}$ and for all $x\in E$, the function $J(MV_{k},V_k)(x,\cdot)$ is continuous on $\mathbb{R}^+$.
\end{prop}
\textbf{Proof} 
We fix $k\in\mathbb{N}$ and $x\in E$. By definition \eqref{def_J} of $J$, we have, for all $t\in\mathbb{R}^+$,
\begin{align*}
J(MV_k,V_k)(x,t)&=
 F(x,t)+ \E_x \left[\right.  \e^{-\alpha (t\wedge t^*(x))} \\
 & MV_k(\Phi(x,t\wedge t^*(x))) \id{S_1 \geq t\wedge t^*(x)}  \\
 &+\e^{-\alpha S_1} V_k(Z_1)\id{S_1 < t\wedge t^*(x)}  \left.\right].
\end{align*}
To show the continuity of $J(MV_{k},V_k)(x,\cdot)$ on $\mathbb{R}^+$, we proceed by the study of the three parts of $J$. The first part of $J$ is the function $F(x,.)$ defined by \eqref{def_F}, which is clearly continuous on $\mathbb{R}^+$. The second part of $J$ is 
\begin{align*}
 \E_x & \left[ \e^{-\alpha (t\wedge t^*(x))}  MV_k (\Phi(x,t\wedge t^*(x))) \id{S_1 \geq t\wedge t^*(x)} \right]  \\
&= \e^{-\alpha (t\wedge t^*(x))}  MV_k (\Phi(x,t\wedge t^*(x)))\mathbb{P}_x( \id{S_1 \geq t\wedge t^*(x)}) \\
&= \e^{-\alpha (t\wedge t^*(x))}  MV_k (\Phi(x,t\wedge t^*(x))) \e^{-\Lambda(x,t\wedge t^*(x))}
\end{align*}
by formula \eqref{loi_T}. 
In the one hand, the function $t\mapsto \e^{-\alpha (t\wedge t^*(x))}\e^{-\Lambda(x,t\wedge t^*(x))}$ is continuous on $\mathbb{R}^+$ by definition \eqref{def_Lambda} of $\Lambda$. In the other hand, using definition \eqref{def_M} of operator $M$, as $\Phi(x,\cdot)$ is continuous on $\mathbb{R}^+$ and $c$ is continuous on $E\times\mathbb{U}$, the function $t\mapsto MV_k(\Phi(x,t\wedge t^*(x)))$ is continuous on $\mathbb{R}^+$. Thus the second part of $J$ is continuous on $\mathbb{R}^+$. \\
The last part of $J$ is continuous on $\mathbb{R}^+$ by the theorem of continuity for a parameter integral. 
Finally, the function $t\mapsto J(MV_{k},V_k)(x,\cdot)$ is continuous on $\mathbb{R}^+$.
\hfill $\Box$


The second result is related to the monotonicity of operators $K$ and $J$ and is useful for the proof of  theorem \ref{ineqVk}.

\begin{lem}\label{ineg_f&g}
Consider $\beta>0$ and suppose that
\begin{equation}\label{ineg01}
\forall x\in E, \qquad \vc{k}(x) \leq V_k(x) \leq \vc{k}(x) + \beta.
\end{equation} 
Then for any $x\in E$ and for any $t\in \mathbb{R}^+$, we have
\begin{equation}\label{ineg-K}
K\vc{k}(x) \leq KV_k(x) \leq K\vc{k}(x) + \beta
\end{equation}
and
\begin{equation}\label{ineg-J}
J(M\vc{k},\vc{k})(x,t) \leq J(MV_k,V_k)(x,t) \leq J(M\vc{k},\vc{k})(x,t) + \beta.
\end{equation}
\end{lem}

\begin{rem}
We have seen that for all $k\in\mathbb{N}$, functions $V_k$ and value functions $ \vc{k}$ are in $\textbf{L}_{\Phi}(E)$. Then functions $K\vc{k}$, $KV_k$, $J(M\vc{k},\vc{k})$ and $ J(MV_k,V_k)$ exist.
\end{rem}

\textbf{Proof}  Those inequalities are showed using the monotonicity of the expectation.
\hfill $\Box$

\section{Markov property of intervention times}

We now prove an important result regarding the Markov property for the intervention times of our strategies.

\begin{prop}\label{prop_markov_tau}
Set $f\in\textbf{B}(\mathbb{R}^+)$ and $\T{x}\in\T{E}$ such that $\T{x}=(x,N,0)$, where $x\in E$ and $N\in\mathbb{N}$. Then for all $i\in\mathbb{N}^*$, we have
\begin{align*}
\T{\E}_{\T{x}}[f(\T{\tau}_i) | \T{\mathcal{F}}_{\T{S}_1}]=&(\id{\T{S}_1<\R{N}(x)\wedge t^*(x)}+\id{\T{S}_1=\R{N}(x)= t^*(x)} \\
&+ \id{\T{S}_1= t^*(x)\neq\R{N}(x)})\T{\E}_{\T{Z}_1}[f(\T{\tau}_i+ \T{S}_1)]  \\
&+\id{\T{S}_1=\R{N}(x)\neq t^*(x)}\T{\E}_{\T{Z}_1}[f(\T{\tau}_{i-1}+ \R{N}(x)].
\end{align*}
\end{prop}

\textbf{Proof} 
\begin{align*}
\T{\E}_{\T{x}}[f(\T{\tau}_i) | &\T{\mathcal{F}}_{\T{S}_1}]\\
=& (\id{\T{S}_1<\R{N}(x)\wedge t^*(x)}+\id{\T{S}_1=\R{N}(x)= t^*(x)} \\
&+ \id{\T{S}_1= t^*(x)\neq\R{N}(x)}+\id{\T{S}_1=\R{N}(x)\neq t^*(x)}) \\
&\times\T{\E}_{\T{x}}\Bigl[\Bigr.f\Bigl( \inf\Big\{t\in\mathbb{R}^+ ;
\sum_{k=1}^{\infty} \id{\T{T}_k\leq t} \id{\T{X}_{\T{T}_k}\in\partial \T{E}}  \\
&\id{\T{S}_k\neq t^*(\pi(\T{Z}_{k-1}))}=i \Big\}\Bigl.\Bigr) \Big| \T{\mathcal{F}}_{\T{S}_1}\Bigr]. 
\end{align*}
We deal with the first term of the sum in the above expectation. First, we can remark that in all cases, as we take $i\neq0$, we have always $\T{\tau}_i\geq \T{S}_1$. Then we obtain $\id{\T{T}_1\leq t}=1$. Now, we distinguish the three cases:
\begin{itemize}
\item If we have $\T{S}_1<\R{N}(x)\wedge t^*(x)$, then $\T{S}_1\neq \T{t}^*(\T{x})$, thus $\T{X}_{\T{T}_1}\notin\partial \T{E}$ and we can conclude that the first term in the sum is zero. 
\item If we have $\T{S}_1=\R{N}(x)= t^*(x)$ or $\T{S}_1= t^*(x)\neq \R{N}(x)$, then $\T{S}_1=t^*(\pi(\T{Z}_0))$, thus the first term in the sum is zero. 
\item If we have $\T{S}_1=\R{N}(x)\neq t^*(x)$, then $\T{S}_1\neq t^*(\pi(\T{Z}_0))$. Furthermore, we have $\T{S}_1=\T{t}^*(\T{x})$ because we have always $\R{N}(x)\leq t^*(x)$. Thus $\T{X}_{\T{T}_1}\in\partial \T{E}$ and we can conclude that the first term in the sum is equal to $1$.
\end{itemize}
So we have
\begin{align*}
\T{\E}_{\T{x}}&[f(\T{\tau}_i) | \T{\mathcal{F}}_{\T{S}_1}]
=(\id{\T{S}_1<\R{N}(x)\wedge t^*(x)}+\id{\T{S}_1=\R{N}(x)= t^*(x)} \\
&+ \id{\T{S}_1= t^*(x)\neq\R{N}(x)}) 
\T{\E}_{\T{x}}\Bigl[\Bigr.f\Bigl(\Bigr. \inf\Big\{t\in\mathbb{R}^+ \Big|
\sum_{k=2}^{\infty} \id{\T{T}_k\leq t} \\
& \id{\T{X}_{\T{T}_k}\in\partial \T{E}} \id{\T{S}_k\neq t^*(\pi(\T{Z}_{k-1}))}=i \Big\}\Bigl.\Bigr) \Big| \T{\mathcal{F}}_{\T{S}_1}\Bigl.\Bigr]  \\
&+\id{\T{S}_1=\R{N}(x)\neq t^*(x)} \T{\E}_{\T{x}}\Bigl[\Bigr.f\Bigl(\Bigr. \inf\Big\{t\in\mathbb{R}^+ \Big|
\sum_{k=2}^{\infty} \id{\T{T}_k\leq t} \\
&\id{\T{X}_{\T{T}_k}\in\partial \T{E}} \id{\T{S}_k\neq t^*(\pi(\T{Z}_{k-1}))}=i-1 \Big\}\Bigl.\Bigr) \Big| \T{\mathcal{F}}_{\R{N}(x)}\Bigl.\Bigr]. 
\end{align*}
Because the process $\{(\T{Z}_n,\T{S}_n)_n\}$ is a Markov chain, the Markov property gives us
\begin{align*}
\T{\E}_{\T{x}}&[f(\T{\tau}_i) | \T{\mathcal{F}}_{\T{S}_1}]
=(\id{\T{S}_1<\R{N}(x)\wedge t^*(x)}+\id{\T{S}_1=\R{N}(x)= t^*(x)} \\
&+ \id{\T{S}_1= t^*(x)\neq\R{N}(x)}) \T{\E}_{\T{x}}\Bigl[\Bigr.f\Bigl(\Bigr. \inf\Big\{t\in\mathbb{R}^+ \Big|
\sum_{k=2}^{\infty} \id{\T{T}_{k-1}\leq t} \\
& \id{\T{X}_{\T{T}_{k-1}}\in\partial \T{E}} \id{\T{S}_{k-1}\neq t^*(\pi(\T{Z}_{k-2}))}=i \Big\} +\T{S}_1\Bigl.\Bigr)\Bigl.\Bigr]  \\
&+\id{\T{S}_1=\R{N}(x)\neq t^*(x)} \T{\E}_{\T{x}}\Bigl[\Bigr.f\Bigl(\Bigr. \inf\Big\{t\in\mathbb{R}^+ \Big|
\sum_{k=2}^{\infty} \id{\T{T}_{k-1}\leq t} \\
&\id{\T{X}_{\T{T}_{k-1}}\in\partial \T{E}} \id{\T{S}_{k-1}\neq t^*(\pi(\T{Z}_{k-2}))}=i-1 \Big\}+\R{N}(x)\Bigl.\Bigr)\Bigl.\Bigr]. 
\end{align*}
Changing the index in the sums and using definition \eqref{def_tps_int} of intervention times, we obtain the result.
\hfill $\Box$

\section{Distributions  of $(\tilde{Z}_1,\tilde{S}_1)$ and $(Z_1,S_1)$}

Finally, we compare the distributions of the first post jump location and sojourn-time for the original PDMP and the auxiliary one.  

\begin{lem}\label{loi_conjointe}
Consider $\T{x}\in\T{E}$ such that $\T{x}=(x,N_0,0)$, where $x\in E$ and $N_0\in\mathbb{N}$. Let $\T{v} \in \textbf{B}(\T{E}\times \mathbb{R}^+)$. We note:
\[\forall x\in E, \forall N \in \mathbb{N}, \forall t\in\mathbb{R}^+, \quad \T{v}(x,N,0,t)=v(x,N,t).\]
Then we have
\begin{align*}
\T{\E}_{\T{x}}&[\T{v}(\T{Z}_1,\T{S}_1)] \\
=& \E_x \Bigl[ \Bigr. \id{N_0=0}v(Z_1,0,S_1) + (\id{S_1<t^*(x)\wedge\R{N_0}(x)} \\
&\times\id{N_0\neq 0} +  \id{S_1=t^*(x)=\R{N_0}(x)}) v(Z_1,N_0-1,S_1)  \\
&+ \id{S_1 \geq \R{N_0}(x) \neq t^*(x)} v(y_\epsilon^{N_0}(x),N_0-1,\R{N_0}(x))  \Bigl.\Bigr].
\end{align*}
\end{lem}
\textbf{Proof} 
The distribution of variables $\T{Z}_1$ and $\T{S}_1$ entails
\begin{align*}
\T{\E}_{\T{x}}&[\T{v}(\T{Z}_1,\T{S}_1)] \\
=&\int_{\T{E}}\int_0^{\T{t}^*(\T{x})} \T{v}(\T{z},\T{s}) \T{\lambda}(\T{\Phi}(\T{x},\T{s})) \e^{-\T{\Lambda}(\T{x},\T{s})} \T{Q}(\T{\Phi}(\T{x},\T{s}),d\T{z}) d\T{s} \\
&+ \int_{\T{E}} \T{v}(\T{z},\T{t}^*(\T{x})) \e^{-\T{\Lambda}(\T{x},\T{t}^*(\T{x}))} \T{Q}(\T{\Phi}(\T{x},\T{t}^*(\T{x})),d\T{z}).
\end{align*}
By construction of the process $\{\T{X}\}$, the random variable $\T{Z}_1$ can be written almost surely by: $\T{Z}_1=(Z,(N_0-1)\wedge 0,0)$, where $Z$ is a random variable on $E$. Using definition \eqref{t_etoile_tilde} of $\tilde{t}^*$, \eqref{phi_tilde} of $\T{\Phi}$, \eqref{lambda_tilde} of $\T{\lambda}$, \eqref{Lambda_tilde} of $\T{\Lambda}$ and \eqref{defqtilde} of $\T{Q}$, we obtain
\begin{align*}
\T{\E}_{\T{x}}[\T{v}(\T{Z}_1,\T{S}_1)]
=& \int_{E}\int_0^{t^*(x)\wedge \R{N_0}(x)} v(z,(N_0-1)\wedge 0,s) \\
& \lambda(\Phi(x,s)) \e^{-\Lambda(x,s)} Q(\Phi(x,s),dz) ds \\
&+ \id{\R{N_0}(x)\geq t^*(x)} \e^{-\Lambda(x,t^*(x))} \\
&\int_{E} v(z,(N_0-1)\wedge 0,t^*(x)) Q(\Phi(x,t^*(x)),dz)  \\
&+ \id{\R{N_0}(x)< t^*(x)} \e^{-\Lambda(x,\R{N_0}(x))} \\
&v(y_\epsilon^{N_0}(x),N_0-1,\R{N_0}(x)).
\end{align*}
Formula \eqref{loi_T} for $t=\R{N_0}$ and the continuity of the law of $S_1$ on $[0;t^*(x))$ give us
\[ \e^{-\Lambda(x,\R{N_0}(x))}\id{\R{N_0}(x)<t^*(x)}= \mathbb{P}_x(S_1\geq \R{N_0}(x)\neq t^*(x)).\]
Furthermore, the law of $S_1$ is such that $\e^{-\Lambda(x,t^*(x))}=\mathbb{P}_x(S_1= t^*(x))$. Hence
\begin{align*}
\T{\E}_{\T{x}}&[\T{v}(\T{Z}_1,\T{S}_1)] \\
=& \E_x \Bigl[ \Bigr. \id{S_1<t^*(x)\wedge\R{N_0}(x)} v(Z_1,(N_0-1)\wedge 0, S_1) \\
&+ \id{\R{N_0}(x)\geq t^*(x)} \id{S_1=t^*(x)} v(Z_1,(N_0-1)\wedge 0,S_1) \\
&+ \id{S_1 \geq \R{N_0}(x) \neq t^*(x)} v(y_\epsilon^{N_0}(x),N_0-1,\R{N_0}(x))  \Bigl.\Bigr].
\end{align*}
By definition \eqref{rk} of $\R{N_0}$, we always have $\R{N_0}(x)\leq t^*(x)$ for $N_0\geq1$ and for $N_0=0$, we decided that $\R{N_0}(x)=+\infty$ in section \ref{sect_def_strategie} . Furthermore, $t^*(x)<+\infty$ because we supposed that $t^*$ is bounded (assumption \ref{hyp2}).Thus we have the equality $\id{\R{N_0}(x)\geq t^*(x)}=\id{\R{N_0}(x)=t^*(x)}+\id{N_0=0}$. The consequence in the expectation is
\begin{align*}
\T{\E}_{\T{x}}&[\T{v}(\T{Z}_1,\T{S}_1)] \\
=& \E_x \Bigl[ \Bigr. \id{S_1<t^*(x)\wedge\R{N_0}(x)} v(Z_1,(N_0-1)\wedge 0, S_1) \\
&+  \id{S_1=t^*(x)=\R{N_0}(x)} v(Z_1,N_0-1,S_1) \\
&+ \id{N_0=0}\id{S_1=t^*(x)} v(Z_1,0,S_1) \\
&+ \id{S_1 \geq \R{N_0}(x) \neq t^*(x)} v(y_\epsilon^{N_0}(x),N_0-1,\R{N_0}(x))  \Bigl.\Bigr] \\
=& \E_x \Bigl[ \Bigr. 
\id{N_0=0}v(Z_1,0,S_1) \\
&+ \id{S_1<t^*(x)\wedge\R{N_0}(x)}\id{N_0\neq 0} v(Z_1,N_0-1, S_1) \\
&+  \id{S_1=t^*(x)=\R{N_0}(x)} v(Z_1,N_0-1,S_1)  \\
&+ \id{S_1 \geq \R{N_0}(x)\neq t^*(x)} v(y_\epsilon^{N_0}(x),N_0-1,\R{N_0}(x))  \Bigl.\Bigr],
\end{align*}
which corresponds to the result.
\hfill $\Box$

\end{document}